\documentclass{amsart}

\usepackage{fancyhdr}
\usepackage{lastpage}
\usepackage{stmaryrd,yhmath}

\newcommand{\bdis}{\begin{displaymath}}
\newcommand{\edis}{\end{displaymath}}
\newcommand{\be}{\begin{equation}}
\newcommand{\ee}{\end{equation}}
\newcommand{\mbb}{\mathbb}
\newcommand{\mcal}{\mathcal}

\newcommand{\vp}{\varphi}

\newcommand{\tg}{\tilde{\gamma}}
\newcommand{\bg}{\bar{\gamma}}

\newtheorem{lemma}[]{Lemma}

\theoremstyle{definition}

\newtheorem{cor}[]{Corollary}

\theoremstyle{remark}
\newtheorem{remark}[]{Remark}

\newtheorem*{mydef1}{{\bf Theorem}}

\numberwithin{equation}{section}

\begin{document}

\title{Jacob's ladders and the nonlocal interaction of the function $|\zeta(1/2+it)|$ with the function $\arg\zeta(1/2+it)$ on the distance
$\sim(1-c)\pi(t)$}

\author{Jan Moser}

\address{Department of Mathematical Analysis and Numerical Mathematics, Comenius University, Mlynska Dolina M105, 842 48 Bratislava, SLOVAKIA}

\email{jan.mozer@fmph.uniba.sk}


\begin{abstract}
In this paper we obtain a new-type formula - \emph{a mixed formula} - which connects the functions $|\zeta(1/2+it)|$ and
$\arg\zeta(1/2+it)$. This formula cannot be obtained in the classical theory of A. Selberg, and, all the less, in the theories of
Balasubramanian, Heath-Brown and Ivic.
\end{abstract}

\maketitle

\section{Introduction}
Let us remind that in the formula
\be \label{1.1}
\zeta\left(\frac{1}{2}+it\right)=\left|\zeta\left(\frac{1}{2}+it\right)\right|e^{i\arg\zeta\left(\frac{1}{2}+it\right)}
\ee
the $\arg\zeta\left(\frac{1}{2}+it\right)$ is defined as follows. If $t=\gamma$, where the $\beta+i\gamma$ is a zero of $\zeta(s)$, the
$\arg\zeta\left(\frac{1}{2}+it\right)$ is obtained by continuous variation along the straight lines joining 2, $2+it$ and $1/2+it$,
starting from the value $\arg\zeta(2)=0$. If $t=\gamma$ then
\bdis
\arg\zeta\left(\frac{1}{2}+i\gamma\right)=\lim_{t\to\gamma^+}\arg\zeta\left(\frac{1}{2}+it\right).
\edis
Let
\be \label{1.2}
S(t)=\frac{1}{\pi}\arg\zeta\left(\frac{1}{2}+it\right), \ S_1(T)=\int_0^T S(t){\rm d}t.
\ee
First of all, there are the asymptotic formulae for the integrals
\bdis
\int_0^T\left|\zeta\left(\frac{1}{2}+it\right)\right|^2{\rm d}t,\ \int_0^T\left|\zeta\left(\frac{1}{2}+it\right)\right|^4{\rm d}t,
\edis
(Hardy and Littlewood started to study these integrals in 1918, 1922). On the other hand, there are Selberg's asymptotic formulae
(\cite{13},\cite{14}, in 1944, 1946)
\be \label{1.3}
\int_T^{T+U}\{ S(t)\}^{2k}{\rm d}t\sim\frac{(2k)!}{k!(2\pi)^{2k}}U(\ln\ln T)^k,
\ee
\be \label{1.4}
\int_T^{T+U}\{ S_1(t)\}^{2k}{\rm d}t\sim c_kU,\ T\to\infty ,
\ee
where $U=T^{1/2+\epsilon}$ and $k$ is a fixed positive number, (Littlewood and Titchmarsh started to study the integrals (\ref{1.3}) and
(\ref{1.4}), $k=1$, see \cite{2},\cite{15} in 1925, 1928). \\

In this paper we obtain a formula of a new type - \emph{a mixed formula} - i.e. a formula which connects the functions (see (\ref{1.1}))
$|\zeta(1/2+it)|, \arg\zeta(1/2+it)$. \\
This paper is a continuation of the series of papers \cite{3}-\cite{12}.

\section{Result}

Let $\vp_1(T),\ T\geq T_0[\vp_1]$ stands for the Jacob's ladder. The following theorem holds true.

\begin{mydef1}
For every fixed $k\in\mbb{N}$ and for every fixed Jacob's ladder there is the single-valued function of $T$
\bdis
\tau_k=\tau_k(T;\vp_1)=\tau_k(T),\ T\geq T_1[\vp_1],
\edis
for which the following asymptotic formula
\be \label{2.1}
\left|\int_0^{\vp_1(\tau_k(T)}\arg\zeta\left(\frac{1}{2}+it\right){\rm d}t\right|\sim \pi(c_k)^{\frac{1}{2k}}
\frac{(\ln\tau_k(T))^{\frac{1}{2k}}}{\left|\zeta\left(\frac{1}{2}+i\tau_k(T)\right)\right|^{\frac{1}{k}}}
\ee
is true, where
\begin{itemize}
\item[(A)] $\tau_k\in (T,T+U),\ \vp_1(\tau_k)\in(\vp_1(T),\ \vp_1(T+U)),\ U=T^{1/2+\epsilon}$
\item[(B)] $\vp_1(T+U)-\vp_1(T)\sim U,\ \vp_1(T+U)< T$
\item[(C)] $\rho\{[T,T+U];[\vp_1(T),\vp_1(T+U)]\}\sim (1-c)\pi(T)\to\infty$ as $T\to\infty$ and $\rho$ denotes the distance of the
corresponding segments, $c$ is the Euler's constant and $\pi(T)$ is the prime-counting function.
\end{itemize}
\end{mydef1}

\begin{remark}
By (\ref{2.1}) we have the prediction of the value
\bdis
\left|\zeta\left(\frac{1}{2}+i\tau_k(T)\right)\right|=|Z(\tau_k(T))|,\ \tau_k(T)\in (T,T+U)
\edis
for the signal
\bdis
Z(t)=e^{i\vartheta(t)}\zeta\left(\frac{1}{2}+it\right)
\edis
generated by the Riemann zeta-function, by means of the value
\bdis
\left|\int_0^{\vp_1(\tau_k(T))}\arg\zeta\left(\frac{1}{2}+it\right){\rm d}t\right|,\ \vp_1(\tau_k(T))\in(\vp_1(T),\vp_1(T+U))
\edis
which descends from very deep past (see (B),(C), comp. \cite{7}, Remarks 3,4).
\end{remark}

\begin{remark}
It is quite evident that the formula (\ref{2.1}) cannot be obtained within the classical theory of A. Selberg (see \cite{13},\cite{14}), and, all
the less, in the theories of Balasubramanian, Heath-Brown and Ivic, (comp. \cite{1}).
\end{remark}

\section{The first corollaries}

Using the mean-value theorem in (\ref{2.1}) and putting $\vp_1(\tau_k(T))\sim\tau_k(T)$, (this follows from $t-\vp_1(t)\sim (1-c)\pi(t)$), we
obtain

\begin{cor}
\be \label{3.1}
|\omega[\arg\zeta]|\sim
\frac{\pi(c_k)^{\frac{1}{2k}}(\ln\tau_k(T))^{\frac{1}{2k}}}{\tau_k(T)\left|\zeta\left(\frac{1}{2}+i\tau_k(T)\right)\right|^{\frac{1}{k}}},
\ee
where
\bdis
\omega[\arg\zeta]=\omega[\arg\zeta(1/2+it);\ t\in [0,\vp_1(\tau_k(T))]]
\edis
denotes the mean-value of $\arg\zeta (1/2+it),\ t\in [0,\vp_1(\tau_k(T))]$ and $\vp_1(\tau_k(T))\in (\vp_1(T),\vp_1(T+U))$.
\end{cor}

\begin{remark}
By (\ref{3.1}) the following holds true: the value $|\zeta(1/2+i\tau_k(T))|$ is asymptotically defined by the mean-value $\omega[\arg\zeta]$, which
descends from the very deep past (comp. Remark 1) and vice-versa.
\end{remark}

\begin{remark}
It is evident that the set of zeroes $t=\gamma$ of the function $\zeta(1/2+it),\ t\in [T,T+U]$ is the exceptional set for the values of the function
$\tau_k=\tau_k(T)$, i.e. $\tau_k(T)\not=\gamma$.
\end{remark}

In connection with this we have the following addition

\begin{cor}
\bdis
\lim_{T\to\infty}|\omega[\arg\zeta]|\cdot |\zeta(1/2+i\tau_k(T))|^{\frac{1}{k}}=0 .
\edis
\end{cor}

\section{The asymptotic formula for the distance $\tg^\prime-\tg$ of some subsequence of consecutive zeroes of $\zeta(1/2+it)$}

\subsection{}

For simplicity, we put $k=1$ in (\ref{3.1})
\be \label{4.1}
|\omega[\arg\zeta]|\cdot \left|\zeta\left(\frac{1}{2}+i\tau_1(T)\right)\right|\sim\pi \sqrt{c_1}\frac{\sqrt{\ln\tau_1(T)}}{\tau_1(T)} .
\ee
Let $t=\tg,\tg^\prime$ denotes the consecutive zeroes of the function $|\zeta(1/2+it)|$, for which
\be \label{4.2}
\tg<\tau_1(T)<\tg^\prime,\ \tau_1(T)\in (T,T+U)
\ee
holds true, and $n(\tg)$ denotes the order of the zero $t=\tg$. Since
\begin{eqnarray} \label{4.3}
& &
\left|\zeta\left(\frac{1}{2}+i\tau_1(T)\right)\right|=|Z[\tau_1(T)]| = \\
& &
=\frac{1}{\{ n(\tg)\}!}\left| Z^{(n(\tg))}[\tau_1^1(T)]\right|\cdot [\tau_1(T)-\tg]^{n(\tg)},\ \tg<\tau_1^1(T)<\tau_1(T) \nonumber
\end{eqnarray}
(see (\ref{4.2})) then from (\ref{4.1}), ($\tau_1(T)\sim\tg$), the asymptotic formula
\be \label{4.4}
\tau_1(T)-\tg\sim
\left(\frac{\pi\sqrt{c_1}\{ n(\tg)\}!\sqrt{\ln\tg}}{\tg|\omega[\arg\zeta]|\cdot |Z^{(n(\tg))}[\tau_1^1(T)]|}\right)^{\frac{1}{n(\tg)}}
\ee
follows, and similarly, we obtain
\begin{eqnarray} \label{4.5}
& &
|\tau_1(T)-\tg^\prime|=\tg^\prime-\tau_1(T)\sim \\
& &
\left(\frac{\pi\sqrt{c_1}\{ n(\tg^\prime)\}!\sqrt{\ln\tg^\prime}}{\tg^\prime|\omega[\arg\zeta]|\cdot |Z^{(n(\tg^\prime))}[\tau_1^1(T)]|}\right)^{\frac{1}{n(\tg^\prime)}}, \
\tau_1(T)<\tau_1^2(T)<\tg^\prime . \nonumber
\end{eqnarray}
Then we obtain from (\ref{4.4}), (\ref{4.5})
\begin{cor}
If the pair of consecutive zeroes $\tg,\tg^\prime$ fulfils (\ref{4.2}) then the asymptotic formula
\begin{eqnarray} \label{4.6}
& &
\tg^\prime-\tg\sim
\left(\frac{\pi\sqrt{c_1}\{ n(\tg)\}!\sqrt{\ln\tg}}{\tg|\omega[\arg\zeta]|\cdot |Z^{(n(\tg))}[\tau_1^1(T)]|}\right)^{\frac{1}{n(\tg)}}+\\
& &
+\left(\frac{\pi\sqrt{c_1}\{ n(\tg^\prime)\}!\sqrt{\ln\tg^\prime}}{\tg^\prime|\omega[\arg\zeta]|\cdot |Z^{(n(\tg^\prime))}[\tau_1^1(T)]|}\right)^{\frac{1}{n(\tg^\prime)}}, \
\tg<\tau_1^1(T)<\tau_1(T)<\tau_1^2(T)<\tg^\prime \nonumber
\end{eqnarray}
holds true.
\end{cor}

\begin{remark}
This is the surprise (for the author) that the asymptotic formula for $\tg^\prime-\tg$ contains the expression $|\omega[\arg\zeta]|$.
\end{remark}

\begin{remark}
The nonlocal nature of the formula (\ref{4.6}) is quite evident (comp. Remark 1).
\end{remark}

\subsection{}

Let $\{ (\bg,\bg^\prime)\}$ denotes the subsequence of the sequence $\{ (\tg,\tg^\prime)\}$ such that $n(\bg)=n(\bg^\prime)=1$
(for example if the weakened Mertens hypothesis is true). Then from (\ref{4.6}) we obtain

\begin{cor}
\begin{eqnarray} \label{4.7}
& &
\bg^\prime-\bg\sim\frac{\pi\sqrt{c_1}\sqrt{\ln\bg}}{\bg|\omega[\arg\zeta]|}
\left(\frac{1}{|Z^\prime[\tau_1^1]|}+\frac{1}{|Z^\prime[\tau_1^2]|}\right),\ \bg<\tau_1^1<\tau_1<\tau_1^2<\bg^\prime .
\end{eqnarray}
\end{cor}

Next, using the estimate
\be \label{4.8}
|Z^\prime|<A\bg^{1/6}\ln\bg,\ \bg^\prime-\bg<A\bg^{1/6+\epsilon/2}
\ee
we obtain from (\ref{4.7}) the following

\begin{cor}
\be \label{4.9}
|\omega[\arg\zeta]|>\frac{A}{\bg^{4/3+\epsilon}},\ \bg\to\infty .
\ee
\end{cor}

\begin{remark}
The estimate (\ref{4.9}) is the first nontrivial lower bound of $|\omega[\arg\zeta]|$, (the small improvements of the exponent 1/6 in
(\ref{4.8}) are irrelevant).
\end{remark}

\section{Lemmas}

\subsection{}

Let us remind that
\be \label{5.1}
\tilde{Z}^2(t)=\frac{{\rm d}\vp_1(t)}{{\rm d}t},\ \vp_1(t)=\frac{1}{2}\vp(t),
\ee
where
\be \label{5.2}
\tilde{Z}^2(t)=\frac{Z^2(t)}{2\Phi^\prime_\vp[\vp(t)]}=\frac{Z^2(t)}{\left\{ 1+\mcal{O}\left(\frac{\ln\ln t}{\ln t}\right)\right\}\ln t},
\ee
(see \cite{3}, (3.9); \cite{5}, (1.3); \cite{9}, (1.1), (3.1), (3.2)). The following Lemma holds true (see \cite{8}, (2.5); \cite{9}, (3.3)).

\begin{lemma}
For every integrable function (in the Lebesgue sense) $f(x),\ x\in [\vp_1(T),\vp_1(T+U)]$ the following is true
\be \label{5.3}
\int_T^{T+U}f[\vp_1(t)]\tilde{Z}^2(t){\rm d}t=\int_{\vp_1(T)}^{\vp_1(T+U)}f(x){\rm d}x,\ U\in\left(\left. 0,\frac{T}{\ln T}\right.\right],
\ee
where $t-\vp_1(t)\sim (1-c)\pi(t)$.
\end{lemma}

\begin{remark}
The formula (\ref{5.3}) remains true also in the case when the integral on the right-hand side of (\ref{5.3}) is only relatively convergent
improper integral of the second kind (in the Riemann sense).
\end{remark}

\subsection{}
Next, the following $\tilde{Z}^2$-transformation of the formulae of A. Selberg (\ref{1.3}), (\ref{1.4}) holds true (comp. \cite{7}, Concluding remarks).

\begin{lemma}
\be \label{5.4}
\int_T^{T+U}\{ S[\vp_1(t)]\}^{2k}\left|\zeta\left(\frac{1}{2}+it\right)\right|^2{\rm d}t\sim\frac{(2k)!}{k!(2\pi)^{2k}}U\ln T(\ln\ln T)^k,
\ee
\be \label{5.5}
\int_T^{T+U}\{ S_1[\vp_1(t)]\}^{2k}\left|\zeta\left(\frac{1}{2}+it\right)\right|^2{\rm d}t\sim c_k U\ln T .
\ee
\end{lemma}
\begin{proof}
From (\ref{1.3}) and by Lemma 1 we obtain
\begin{eqnarray} \label{5.6}
& &
\int_T^{T+U}\{ S[\vp_1(t)]\}^{2k}\tilde{Z}^2(t){\rm d}t=\int_{\vp_1(T)}^{\vp_1(T+U)}\{ S[\vp_1(T)]\}^{2k}{\rm d}t\sim \\
& &
\sim \frac{(2k)!}{k!(2\pi)^{2k}}\{ \vp_1(T+U)-\vp_1(T)\}\ln\ln \vp_1(T)= \nonumber \\
& &
=\frac{(2k)!}{k!(2\pi)^{2k}} U\frac{\vp_1(T+U)-\vp_1(T)}{U}\ln\ln \vp_1(T)= \nonumber \\
& &
=\frac{(2k)!}{k!(2\pi)^{2k}} U \tan[\alpha(T,U)]\ln\ln \vp_1(T) , \nonumber
\end{eqnarray}
where $\alpha(T,U)$ is the angle of the chord of the curve $y=\vp_1(t)$ that binds the points $[T,\vp_1(T)]$ and $[T+U,\vp_1(T+U)]$. Let us remind that
\be \label{5.7}
\tan[\alpha(T,U)]\sim 1,\ U\in \left[ T^{1/3+\epsilon},\frac{T}{\ln T}\right],
\ee
(comp. \cite{3}, (8.3), \cite{7}, (3.9)) and
\be \label{5.8}
t-\vp_1(t)\sim (1-c)\pi(t) ,
\ee
(see (5.1), \cite{3}, (6.2)). Since by (5.8) we have
\be \label{5.9}
\vp_1(T)\sim T,\ \ln\ln \vp_1(T)\sim \ln\ln T,
\ee
then from (5.6) by (5.7), (5.9) this formula
\be \label{5.10}
\int_T^{T+U}\{ S[\vp_1(t)]\}^{2k}\tilde{Z}^2(t){\rm d}t\sim\frac{(2k)!}{k!(2\pi)^{2k}} U (\ln\ln T)^k
\ee
follows. Using the mean-value theorem in (\ref{5.10}) we obtain (\ref{5.4}) by (\ref{5.2}). Similarly we obtain (\ref{5.5}).
\end{proof}

\section{Proof of the Theorem}

Using the mean-value theorem in (\ref{5.5}) we obtain
\begin{eqnarray} \label{6.1}
& &
|S_1[\vp_1(\tau)]|^{2k}Z^2(\tau)\sim c_k\ln T\sim c_k\ln\tau , \\
& &
\tau=\tau(T,U,k)\in (T,T+U),\ \vp_1(\tau)\in (\vp_1(T),\vp_1(T+U)) . \nonumber
\end{eqnarray}
Since $U=T^{1/2+\epsilon}$ then $\tau=\tau(T,k,\epsilon),\ T\geq T_0[\vp_1]$. Next, if for every $T\geq T_0[\vp_1]$ we take one mean-value
$\tau_k(T)\in \{ \tau(T,k,\epsilon)\}$ then we have the single-valued function $\tau_k=\tau_k(T)$ for fixed $k,\epsilon$
(comp. \cite{7}, Remark 2). Hence by (\ref{5.1}) we have
\begin{eqnarray} \label{6.2}
& &
|S_1[\vp_1(\tau_k(T))]|^{2k}\sim c_k\frac{\ln \tau_k(T)}{\left|\zeta\left(\frac{1}{2}+i\tau_k(T)\right)\right|^2},\ \\
& &
\tau_k(T)\in (T,T+U),\ \vp_1(\tau_k(T))\in (\vp_1(T),\vp_1(T+U)). \nonumber
\end{eqnarray}
Then from (\ref{6.2}) the asymptotic formula (\ref{2.1}) and (A) follows. The expressions (B) and (C) are identical with (C), (D), $k=1$ of
Theorem in \cite{9}.

\thanks{I would like to thank Michal Demetrian for helping me with the electronic version of this work.}

\end{document}